 \documentclass[final,1p,times]{elsarticle}
\usepackage{listings}
\usepackage{color}
\usepackage[]{multicol}
\usepackage{amsthm}
\usepackage{caption}
\usepackage{here}
\usepackage{mathrsfs}
\usepackage{subcaption}
\usepackage{amsmath, amssymb}
\usepackage{ascmac}
\usepackage{url}
\usepackage{comment}
\usepackage{nccmath}
\usepackage {booktabs}

\newtheorem{theorem}{Theorem}
 
\newtheorem{definition}{Definition}
\newtheorem{lemma}{Lemma}
 
\newtheorem{corollary}{Corollary}

\begin{document}

\begin{frontmatter}



\title{Numerical computation for the exact distribution of \\Roy's largest root statistic under linear alternative}


\author[label1]{Koki Shimizu}
\author[label1]{Hiroki Hashiguchi}

\address[label1]{Tokyo University of Science, 1-3 Kagurazaka, Shinjuku-ku, Tokyo, 162-8601, Japan}
\cortext[mycorrespondingauthor]{Corresponding author. Email address: \url{1420702@ed.tus.ac.jp}~(K. Shimizu).}

\begin{abstract}
This paper discusses the computation of exact powers for Roy's test in multivariate analysis of variance~(MANOVA). 
We derive an exact expression for the largest eigenvalue of a singular noncentral Beta matrix in terms of the product of zonal polynomials.
The numerical computation for that distribution is conducted by an algorithm that expands the product of zonal polynomials as a linear combination of zonal polynomials.
Furthermore, we provide an exact distribution of the largest eigenvalue in a form that is convenient for numerical calculations under the linear alternative. 

\end{abstract}




\end{frontmatter}


\section{Introduction}
\noindent
The distributions of central and noncentral matrix variate beta type I and II play an important role in multivariate analysis.
These distributions have been well-studied in the case of a nonsingular matrix, see Muirhead (1982), Srivastava and Khatri~(1979), and Gupta and Nager~(2000).
The matrix variate beta type I and II distributions were extended to the case of a singular matrix by D{\'\i}az-Garc{\'\i}a and Guti{\'e}rrez-J{\'a}imez~(1997), Srivastava~(2003) and D{\'\i}az-Garc{\'\i}a and Guti{\'e}rrez-J{\'a}imez~(2008). 
Matrix variate beta type II is also known as matrix variate $F$-distribution.
Matrices that follow the beta-type I and II distributions are called Beta and $F$ matrices, respectively.
Wilk's Lambda, Lawley-Hotelling trace, Bartlett-Nanda-Pillai trace, and Roy's largest root are the statistics for testing problems in the MANOVA, canonical correlation analysis, general linear model, and test of equality of covariance matrices.
The four statistics are functions of eigenvalues for the central or noncentral Beta and $F$ matrices.

Roy's largest root is the only statistic that uses one eigenvalue. Therefore, Roy's test is more powerful than the other tests under the linear alternative; i.e., the groups are heavily concentrated in a line. 
Under the null hypothesis in the MANOVA problem, there are efficient calculation and approximation methods for the distribution of Roy's largest root statistic.
Chiani~(2016) determined the exact expression as a Pfaffian of a skew-symmetric matrix and provided an algorithm to compute it.
Johnstone~(2008) showed that Tracy-Widom distribution approximates the exact expression.
These distributions are derived from the distribution of eigenvalues of a nonsingular central Beta matrix.
Ang et al.~(2021) compared Chiani's algorithm with Johnstone's Tracy-Widom approximation to determine the conditions under which each method is appropriate. 

In contrast, the exact computation for the distribution of Roy's largest root statistic under the alternatives has long been an open problem.
Pillai and Sugiyama (1969) derived the exact distribution of the largest eigenvalue for a nonsingular noncentral Beta matrix in terms of the product of the zonal polynomials. 
However, the computation for that distribution could not be performed because the algorithm for computing the product of zonal polynomials remained unknown.
Recently, Shimizu and Hashiguchi~(2022b) provided that algorithm and applied it to compute the distribution of extreme eigenvalues of a singular Wishart matrix for the sphericity test. 
Under the linear alternative, Johnstone and Nader (2017) presented the approximate distribution of Roy's largest root statistic by a combination of independent central and noncentral $F$ variates.
However, this approximation is not accurate when the dimension is large and the noncentrality parameter is small.

In this study, we compute the exact distribution of Roy's largest root statistic under the alternatives.
Section~\ref{sec:02} discusses the problem of expressing the product of two zonal polynomials by a linear combination. 
In Section~\ref{sec:03}, we derive the exact distribution of the largest eigenvalue of a singular noncentral Beta matrix in terms of the product of zonal polynomials.
This is an extension of the results found by Pillai and Sugiyama (1969) for a nonsingular noncentral Beta matrix.
Furthermore, we present the distribution in a form that does not require the computation of zonal polynomials under the linear alternative.
Finally, we compute the theoretical distribution and obtain the exact powers for Roy's test in Section~\ref{sec:04}. 

\section{Preliminary}
\label{sec:02}
In this section, we discuss the calculation that expresses the product of zonal polynomials as a linear combination. 
For a positive integer $k$, let $\kappa=(\kappa_1,\kappa_2,\dots,\kappa_l)$ denote a partition of $k$ with $\kappa_1\geq\cdots\geq \kappa_l> 0$ and $\kappa_1+\cdots +\kappa_l=k$, where $l=l(\kappa)$ is the partition length.
It is often convenient to consider $\kappa$ as having any number of additional zeros $\kappa=(\kappa_1,\kappa_2,\dots,\kappa_l, 0,\dots, 0)$.
The set of all partitions with lengths less than or equal to $m$ is denoted by $P^k_{m}=\{\ \kappa=(\kappa_1,\dots,\kappa_m)\mid \kappa_1+\dots+\kappa_m=k, \kappa_1\geq \kappa_2\geq\cdots\geq  \kappa_m \geq 0 \}$. 
Furthermore, let $X$ be an $m\times m$ symmetric matrix with eigenvalues $x_1,\dots,x_m$. 
 The elementary symmetric functions in eigenvalues $x_1, \dots, x_m$ of $X$ are expressed as $e_0=1$ and 
\begin{align*}
e_1=x_1+\cdots +x_m, e_2=x_1x_2+\cdots+x_{m-1}x_m, \ldots, e_m=x_1x_2\cdots x_m.
\end{align*}
For $\kappa\in P^k_{m}$, we define the polynomials $\mathcal{E}_\kappa(X)$ as
\begin{align*}
\mathcal{E}_\kappa(X)=e_1^{\kappa_1-\kappa_2}e_2^{\kappa_2-\kappa_3}\cdots e_{m-1}^{\kappa_{m-1}-\kappa_{m}}e_m^{\kappa_m}, 
\end{align*}
where the degree of $\mathcal{E}_\kappa(X)$ is $(\kappa_1-\kappa_2)+2(\kappa_2-\kappa_3)+\cdots+m\kappa_m=k$.
The definition of the polynomials $\mathcal{E}_\kappa(X)$ is given in Hashiguchi et al.~(2000) and Jiu and Koutschan~(2020). 
In Muirhead (1982) and Dumitriu et al.~(2007), the zonal polynomial was defined as an expansion of the monomial symmetric polynomial.
In contrast, Hashiguchi et al.~(2000) defined the zonal polynomial in terms of $\mathcal{E}_\mu$'s as follows.
\begin{definition}\label{def1} 
 For an $m\times m$ symmetric matrix $X$ and $\kappa \in P^k_m$, there exists a unique polynomial $C_\kappa(X)$ that satisfies the following three conditions: 
\begin{itemize}
\item[1] $C_\kappa(X)=\sum_{\mu\preceq \kappa} q[\kappa, \mu]\mathcal{E}_\mu(X)$ and $q[\kappa, \kappa] \neq 0$,
where $q[\kappa, \mu]$ is a constant, and the summation is over all partitions $\mu\preceq \kappa$; that is, $\mu$ is less than or equal to $\kappa$ in the lexicographical ordering. 

 \item[2] $D_{m}C_\kappa(X)=d(\kappa) C_\kappa(X)$,
 where $D_{m}$ is the Laplace--Beltrami operator
 $ D_{m}=
    \sum_{i=1}^{m}x_{i}^{2}\frac{\partial^{2}}{\partial x_{i}^{2}}
    +\sum_{i\leq i \neq j \leq m}\frac{x_{i}^{2}}{x_{i}-x_{j}}\frac{\partial}{\partial x_{i}}$, $d(\kappa)$ is the corresponding eigenvalue
    $d(\kappa)=\sum_{i=1}^{m}\kappa_i(\kappa_i+m-i-1)$.
    
 \item[3] $[\mathrm{tr}(X)]^k=\sum_{\kappa\in P^k_m}C_\kappa (X)$.
\end{itemize}
\end{definition}
 We refer to the above polynomials $C_\kappa(X)$ as zonal polynomials. 
 Hashiguchi et al.~(2000) proposed an algorithm that expands the zonal polynomials in terms of $\mathcal{E}_\mu$'s.
If $X=I_m$, the zonal polynomials $C_\kappa(X)$ can be represented by 
\begin{align}
\label{zonal-identity}
C_{\kappa}(I_{m})=2^{2k}k!(m/2)_\kappa\frac{\prod_{i<j}^{l(\kappa)}(2\kappa_i-2\kappa_j-i+j)}{\prod_{i=1}^{l(\kappa)}(2\kappa_i+l(\kappa)-1)!},
\end{align}
where $(\alpha)_\kappa=\prod_{i=1}^{l(\kappa)}\{\alpha-(i-1)/2\}_{\kappa_i}$, and $(\alpha)_k=\alpha(\alpha+1)\cdots (\alpha+k-1)$, and $(\alpha)_0=1$. 

For $\kappa\in P^k_m$ and $\tau\in P^t_m$, we express the product of zonal polynomials as a linear combination of such polynomials:
\begin{align}
\label{prod-zonal}
C_\kappa(X)\cdot C_\tau(X)=\sum_{\delta\in P^{k+t}_m} g^\delta_{\kappa,\tau}C_\delta(X),
\end{align}
where $g^\delta_{\kappa,\tau}$ is a constant determined by the partitions $\kappa$ and $\tau$.
The right-hand side of \eqref{prod-zonal} appears in certain areas of multivariate statistics, such as the sphericity test, instrumental variable~(IV) regression, and power functions in MANVOA; see Sugiyama~(1970), Shimizu and Hashiguchi~(2022b), Hillier and Kan~(2021) and Pillai and Sugiyama~(1969).
The coefficients $g^\delta_{\kappa,\tau}$ were tabulated for values up to the $7$th degree of $k$ and $t$ in Hayakawa~(1967) and Khatri and Pillai~(1968).
Shimizu and Hashiguchi (2022b) proposed an algorithm for obtaining the coefficients $g^\delta_{\kappa,\tau}$ in \eqref{prod-zonal}. 
The symbolic computation of the product $\mathcal{E}_\kappa\cdot\mathcal{E}_\tau$  is simple to carry out.
Thus, the important aspect is to express the zonal polynomials in terms of $\mathcal{E}_\mu$'s.
If $\kappa=(2)$, $\tau=(3,2)$, and $m=2$, the product of the zonal polynomials is represented by 
\begin{align*}
C_{(2)}(X)\cdot C_{(3,2)}(X)&=\biggl(\mathcal{E}_{(2)}(X)-\frac{4}{3}\mathcal{E}_{(1,1)}(X)\biggl) ~\cdot \frac{48}{7}\mathcal{E}_{(3,2)}(X)\\\nonumber
&=\frac{48}{7}\mathcal{E}_{(5,2)}(X)-\frac{64}{7}\mathcal{E}_{(4,3)}(X)\\ 
&=\frac{99}{245}C_{(5,2)}(X)+\frac{12}{35}C_{(4,3)}(X)
\end{align*}
Here, the third line is satisfied by the transition algorithm from the polynomials $\mathcal{E}_{\kappa}(X)$ to the linear combination of $\{C_{\kappa}(X) \mid \kappa \in P_{m}^{k + t}\}$ obtained by Shimizu and Hashiguchi~(2022b).
Another approach for calculating the product of zonal polynomials was presented by Hillier and Kan~(2021).
Here, the zonal polynomials were expressed in terms of power sum functions, and the desired result was obtained using the transition matrix between power sum functions and zonal polynomials.
The calculation of the product of the power sum functions is also straightforward, as in the case of elementary symmetric functions.
However, it requires a considerable amount of calculation time to expand the zonal polynomials to large degrees.
The following lemma is a formula for the special case of $g^\delta_{\kappa,\tau}$ given by Kushner (1988).
\begin{lemma}
Let $\delta=(\delta_1,\dots, \delta_l)$ and $\tau=(\tau_1,\dots,\tau_l)$ be the partitions of $k+t$ and $t$, respectively.
If $\delta_i\geq \tau_i \geq \delta_{i+1}$, $1\leq i\leq l(\tau)$,
then the coefficients $g^{\delta}_{{(k)},\tau}$ are represented by
 \begin{align}
 \label{g-coeffcient}
 \nonumber
g^{\delta}_{{(k)},\tau}=&\frac{k!}{(2k)! \binom{k+t}{k}}\prod_{i=1}^{l(\tau)+1}\frac{[2(\delta_i-\tau_i)-1,1]_2}{(\delta_i-\tau_i)!}\prod_{i<j}^{l(\tau)+1}(T_i-T_j)\\
&\times \prod_{i<j}^{l(\tau)+2}[D_i-T_j-1, T_i-D_j+1]_2\bigg/ \prod_{i<j}^{l(\tau)+1}[D_i-T_j, T_i-D_j]_2,
\end{align}
where $T_i=2\tau_i-i$ and $D_i=2\delta_i-i$.
The factorial symbol $[b, c]_2$ is defined by 
\begin{align*}
[b, c]_2&=b(b-2)(b-4)\cdots c ~~ \text{if}~b-c\geq 2n\\
[b, c]_2&=1~~~~~~~~~~~~~~~~~~~~~~~~~~~~~~~ \text{if}~b-c=-2 
\end{align*}
\end{lemma}
The formula of $g^{\delta}_{{(k)},\tau}$ is useful when computing the exact distribution of the Bartlett-Nanda-Pillai trace and Roy's largest root statistics under the linear alternative.
\section{Exact distribution of the largest eigenvalue under alternative in MANOVA}
\label{sec:03}
We consider the test of equality for mean vectors from $p$ groups.  
Let $x_{ij}$ be distributed as the $i$-th population $N_m(\mu_i, \Sigma)$, $i=1,\dots, p$, $j=1,\dots, n_i$, where $n-p\geq m$ and $n=\sum_{i=1}^{p}n_i$.
We assume that all covariances are equal.
We test the equality of mean vectors; $H_{0}: \mu_1=\cdots =\mu_p$, versus the alternative hypnosis $H_1$ that $\mu_i$ are not all equal. 
The between-group and within-group covariance matrices are defined by
\begin{align*}
H=\sum_{i=1}^{p}n_i(\bar{x}_i-\bar{x})(\bar{x}_i-\bar{x})^\top, ~~~E=\sum_{i=1}^{p}\sum_{j=1}^{n_i}(x_{ij}-\bar{x}_i)(x_{ij}-\bar{x}_i)^\top,
\end{align*}
respectively, where $\bar{x}_i=n^{-1}_i\sum_{j=1}^{n_i}x_{ij}$ and $\bar{x}=\sum_{i=1}^{p}\sum_{j=1}^{n_i}x_{ij}/n$.
The matrix $E$ has a nonsingular Wishart distribution with $n_E$ degrees of freedom, denoted by $E\sim W_m(n_E, \Sigma)$, where $n_E=n-p$.
Under $H_1$, the matrix $H$ has a nonsingular or singular non-central Wishart distribution, denoted by $H\sim W_m(n_H,\Sigma,\Omega)$, where $n_H=p-1$, $\Omega=\Sigma^{-1}\sum_{i=1}^{p}n_i(\mu_i-\bar{\mu})(\mu_i-\bar{\mu})^\top$ and $\bar{\mu}=(1/n)\sum_{i=1}^{p}\mu_i$. (See Fujikoshi et al.~(2010), pp.157 and Muirhead~(1982), pp.486).
Throughout this paper, the notation $W_m(n_H,\Sigma,\Omega)$ is referred to in both nonsingular and singular cases.
Let $H\sim W_m(n_H,\Sigma,\Omega)$ and $E\sim W_m(n_E, \Sigma)$, where $H$ and $E$ are independent.
The noncentral Beta matrix $U$ can be defined as $U=(H+E)^{-1/2}H(H+E)^{-1/2}$.
If $m>n_H$, the matrix $U$ is singular; otherwise, it is nonsingular.
The spectral decomposition of $U$ is $U=H_1LH_1^\top$, where $L=\mathrm{diag}(\ell_1,\dots, \ell_{n_\mathrm{min}})$, $1>\ell_1>\cdots >\ell_{n_\mathrm{min}}>0$, $n_\mathrm{min}=\mathrm{min}\{n_H, m\}$ and $V_{n_\mathrm{min},m}=\{H_1\mid H_1^\top H_1=I_{n_\mathrm{min}}\}$.
The largest eigenvalue $\ell_1$ and sum of eigenvalues $V^{(n_\mathrm{min})}=\sum_{i=1}^{n_\mathrm{min}}\ell_i$ are called the Roy's largest root and Bartlett-Nanda-Pillai trace statistics, respectively.
D{\'\i}az-Garc{\'\i}a and Guti{\'e}rrez-J{\'a}imez~(2008) gave the density function of $U$ with $n_E \geq m>n_H$ as
\begin{align*}
f(U)=C~|L|^{(n_H-m-1)/2}|I_m-U|^{(n_E-m-1)/2}\mathrm{etr}\biggl(-\frac{1}{2}\Omega\biggl)\sum_{k=0}^{\infty}\sum_{\kappa \in P^m_{n_H}}\frac{\{(n_H+n_E)/2\}_\kappa}{(n_H/2)_\kappa (n_E/2)_\kappa } \frac{C_\kappa(\Omega/2)C_\kappa(U)}{C_\kappa(I_m)},
\end{align*}
where $C=\frac{\pi^{(-mn_H+n_H^2)/2}\Gamma_m(n_H+n_E/2)}{\Gamma_{n_H}(n_H/2)\Gamma_m(n_E/2)}$, $\mathrm{etr}(\cdot)=\mathrm{exp}(\mathrm{tr}(\cdot))$ and $\Gamma_m(a)=\pi^{m(m-1)/4}\prod_{i=1}^{m}\Gamma(a-(i-1)/2)$.
The singular non-central matrix variate beta type II distribution was also given in D{\'\i}az-Garc{\'\i}a and Guti{\'e}rrez-J{\'a}imez~(2008). 
For integers $n_E \geq m>n_H$. 
D{\'\i}az-Garc{\'\i}a and Guti{\'e}rrez-J{\'a}imez~(2008) also gave the joint density of eigenvalues in a singular matrix $U$ as
\begin{align}
\label{jointdensity}
\nonumber
f(\ell_1,\dots, \ell_{n_H})=&C_1~\mathrm{etr}\biggl(-\frac{1}{2}\Omega\biggl)|L|^{(m-n_H-1)/2}|I_{n_H}-L|^{(n_E-m-1)/2}\\
&\times \prod_{i<j}^{}(\ell_i-\ell_j)\sum_{k=0}^{\infty}\sum_{\kappa\in P^k_{n_H}}\frac{\{(n_H+n_E)/2\}_\kappa}{(n_H/2)_\kappa}\frac{C_\kappa(\Omega/2)C_\kappa(L)}{C_\kappa(I_m)k!},
\end{align}
where $C_1=\frac{\pi^{n_H^2/2}\Gamma_m\{(n_H+n_E)/2\}}{\Gamma_{n_H}(n_H/2)\Gamma_{m}(n_E/2)\Gamma_{n_H}(m/2)}$.
The following lemma given by Sugiyama~(1967) is necessary to integrate \eqref{jointdensity} with respect to $\ell_2,\dots, \ell_{n_H}$.
\begin{lemma}
\label{Sugiyama-formula}
\begin{align} 
&\int_{1>x_2>\cdots>x_m>0}|X_2|^{t-(m+1)/2}C_\kappa(X_1)\prod_{i=2}^{m}(1-x_i)\prod_{i<j}(x_i-x_j)\prod_{i=2}^m dx_i\nonumber\\
=&(mt+k)
\frac{\Gamma_m(m/2)(t)_\kappa\Gamma_m(t)\Gamma_m\{(m+1)/2\}}
{\pi^{m^2/2}\{t+(m+1)/2\}_\kappa\Gamma_m\{t+(m+1)/2\}}C_\kappa (I_m),
\end{align}
where $ \mathrm{Re}(t) >(m-1)/2$.
\end{lemma}
The following theorem is the extension of the result of Pillai and Sugiyama (1969) to the singular case.
\begin{theorem}
\label{theo01}
Let $H\sim W_m(n_H,\Sigma,\Omega)$ and $E\sim W_m(n_E, \Sigma)$, where $m>n_H$ and $H$ and $E$ independent;
then, the distribution function of the largest eigenvalue $\ell_1$ of $U=(H+E)^{-1/2}H(H+E)^{-1/2}$ is given as 
\begin{align}
\label{prob-ell1}
\nonumber
\mathrm{Pr}(\ell_1 < x) =& \frac{\Gamma_{n_H}\{(n_H+1)/2\} \Gamma_m\{(n_H+n_E)/2\}}{\Gamma_m(n_E/2) \Gamma_{n_H}\{(n_H + m + 1)/2\}}\mathrm{etr}\biggl(-\frac{1}{2}\Omega\biggl)\sum_{k=0}^{\infty}\sum_{\kappa\in P^k_{n_H}}\frac{\{(n_H+n_E)/2\}_\kappa}{(n_H/2)_\kappa}\frac{C_\kappa(\Omega/2)}{C_\kappa(I_m)k!}  \\
&\times \sum_{t=0}^{\infty} \sum_{\tau\in P^t_{n_H}}\sum_{\delta\in P^{k+t}_{n_H}}g_{\kappa,\tau}^{\delta}\frac{\{(m+1-n_E)/2\}_\tau(m/2)_{\delta}}{\{(n_H+m+1)/2\}_\delta t!}C_\delta (I_{n_H})\cdot x^{mn_H/2+k+t}.
\end{align}
\begin{proof}
This proof is similar to the method presented by Pillai and Sugiyama~(1969).
From \eqref{prod-zonal}, we write
\begin{align*}
|I_{n_H}-L|^{(n_E-m-1)/2} C_\kappa(L)&={_1F_0}\biggl(\frac{m-n_E+1}{2},L\biggl)C_\kappa(L)\\
&=\sum_{t=0}^{\infty}\sum_{\tau\in P^t_{n_H}} \frac{\{(m-n_E+1)/2\}_\tau C_\tau (L)}{t!}C_\kappa(L)\\
&=\sum_{t=0}^{\infty}\sum_{\tau\in P^t_{n_H}}\sum_{\delta \in P^{k+t}_{n_H}}\frac{\{(m-n_E+1)/2\}_\tau}{t!}g^\delta_{\kappa,\tau}C_\delta(L). 
\end{align*}
Therefore, the density function \eqref{jointdensity} is represented by 
\begin{align*}
f(\ell_1,\dots, \ell_{n_H})=&~C_1~\mathrm{etr}\biggl(-\frac{1}{2}\Omega\biggl)|L|^{(m-n_H-1)/2}
\sum_{k=0}^{\infty}\sum_{\kappa\in P^k_{n_H}}\frac{\{(n_H+n_E)/2\}_\kappa}{(n_H/2)_\kappa}\frac{C_\kappa(\Omega/2)}{C_\kappa(I_m)k!}\\
&\times \prod_{i<j}^{n_H}(\ell_i-\ell_j)\sum_{t=0}^{\infty}\sum_{\tau\in P^t_{n_H}}\sum_{\delta \in P^{k+t}_{n_H}}g^\delta_{\kappa,\tau} \frac{\{(m+1-n_E)/2\}_\tau C_\delta(L)}{t!}.
\end{align*}
Translating $x_i=\ell_i/\ell_1$, $i=2,\dots,\ell_{n_H}$ and using Lemma~\ref{Sugiyama-formula}, we have
\begin{align*}
f(\ell_1)=&~C_1~\mathrm{etr}\biggl(-\frac{1}{2}\Omega\biggl)\sum_{k=0}^{\infty}\sum_{\kappa\in P^k_{n_H}}\frac{\{(n_H+n_E)/2\}_\kappa}{(n_H/2)_\kappa}\frac{C_\kappa(\Omega/2)}{C_\kappa(I_m)k!}\sum_{t=0}^{\infty}\sum_{\tau\in P^t_{n_H}}\sum_{\delta \in P^{k+t}_{n_H}}g^\delta_{\kappa,\tau} \frac{\{(m+1-n_E)/2\}_\tau }{t!}\\
&\times \ell_1^{m{n_H}/2+k+t-1}  \int_{1>x_2>\cdots>x_{n_H}>0}|X_2|^{(m-n_H-1)/2}\prod_{i=2}^{n_H}(1-x_i)\prod_{2\leq i<j}^{n_H}(\ell_i-\ell_j)C_\delta(X_1)\prod_{i=2}^{n_H}dx_i\\
=&~\frac{\Gamma_m\{(n_H+n_E)/2\}\Gamma_{n_H}\{({n_H}+1)/2\}}{\Gamma_{m}(n_E/2)\Gamma_{n_H}\{(m+{n_H}+1)/2\}}~\mathrm{etr}\biggl(-\frac{1}{2}\Omega\biggl)\sum_{k=0}^{\infty}\sum_{\kappa\in P^k_{n_H}}\frac{\{(n_H+n_E)/2\}_\kappa}{(n_H/2)_\kappa}\frac{C_\kappa(\Omega/2)}{C_\kappa(I_m)k!}\\
& \times \sum_{t=0}^{\infty}\sum_{\tau\in P^t_{n_H}}\sum_{\delta \in P^{k+t}_{n_H}}g^\delta_{\kappa,\tau} \frac{\{(m+1-n_E)/2\}_\tau }{t!}\ell_1^{m{n_H}/2+k+t-1} (m{n_H}/2+k+t)\frac{(m/2)_\delta \mathcal{C}_\delta (I_{n_H})}
{\{(m+{n_H}+1)/2\}_\delta}
\end{align*}
Finally, integrating $f(\ell_1)$ with respect to $\ell_1$, we have the desired result. 
\end{proof}
\end{theorem}
We can also consider the exact distribution of Roy's largest root statistic for $m>n_H$ from the distribution given by Pillai and Sugiyama~(1969) by replacing $m\to n_H$, $n_H\to m$, and $n_E\to n_H+n_E-m$. 
We note that the distribution is different from \eqref{prob-ell1}.
The generalization for the distribution of $\ell_1$ in $U$ is given as 
\begin{align*}
\mathrm{Pr}(\ell_1 < x) =& \frac{\Gamma_{n_\mathrm{min}}\{(n_\mathrm{min}+1)/2\} \Gamma_m\{(n_H+n_E)/2\}}{\Gamma_m(n_E/2) \Gamma_{n_\mathrm{min}}\{(n_H + m + 1)/2\}}\mathrm{etr}\biggl(-\frac{1}{2}\Omega\biggl)\sum_{k=0}^{\infty}\sum_{\kappa\in P^k_{n_\mathrm{min}}}\frac{\{(n_H+n_E)/2\}_\kappa}{(n_H/2)_\kappa}\frac{C_\kappa(\Omega/2)}{C_\kappa(I_m)k!}  \\
&\times \sum_{t=0}^{\infty} \sum_{\tau\in P^t_{n_\mathrm{min}}}\sum_{\delta\in P^{k+t}_{n_\mathrm{min}}}g_{\kappa,\tau}^{\delta}\frac{\{(m+1-n_E)/2\}_\tau(n_\mathrm{max}/2)_{\delta}}{\{(n_H+m+1)/2\}_\delta t!}C_\delta (I_{n_\mathrm{min}})\cdot x^{mn_H/2+k+t},
\end{align*}
where $n_\mathrm{min}=\{n_{n_H}, m\}$ and $n_\mathrm{max}=\{n_{n_H}, m\}$.
If $n\geq m$, the above distribution coincides with the results found by Pillai and Sugiyama~(1969).

Next, we consider the exact distribution for Roy's test under particular alternatives.
Rencher and Christensen~(2012) proposed two alternative hypotheses as follows.
\begin{enumerate}
\item[{\bf a}] {\bf Linear alternative hypothesis:}
If the largest eigenvalue $\ell_1$ is larger than the others, all mean vectors lie close to a line. 
All mean vectors $\mu_i$ for $i=1,\dots, p$ are expressed with $\mu_i=\xi\mu_0$, where $\mu_0$ is an unknown mean vector and $\xi$ is a scale factor.
In this case, the rank of the non-central matrix is $\mathrm{rank}(\Omega)=1$.
\item[{\bf b}]{\bf Planar alternative hypothesis:}
If the eigenvalues $\ell_1$ and $\ell_2$ are larger, the mean vectors mostly lie in a plane.
All mean vectors $\mu_i$ for $i=1,\dots, p$ are expressed as $\mu_i=\alpha_i\mu_1+\beta_i\mu_2$, where $\mu_1$ and $\mu_2$ are unknown mean vectors and $\alpha_i$ and $\beta_i$ are scale factors.
In this case, the rank of the non-central matrix is $\mathrm{rank}(\Omega)=2$.
\end{enumerate} 
Johnstone and Nadler (2017) approximated the distribution of $\ell_1$ under the linear alternative. 
The exact distributions of Wilks's statistics under the linear and planer alternatives were considered in Phong et al.~(2019).
Corollary~\ref{exact-rank1kei} provides the exact distribution \eqref{prob-ell1} under the linear alternative.
\begin{corollary}
\label{exact-rank1kei}
If $\mathrm{rank}(\Omega)=1$ and $\theta_1$ is the largest eigenvalue of $\Omega$ in Theorem~\ref{theo01}, then we have  
\begin{align}
\label{exact-rank1}
\nonumber
\mathrm{Pr}(\ell_1 < x) 
=& \frac{\Gamma_{n_H}\{(n_H+1)/2\} \Gamma_m\{(n_H+n_E)/2\}}{\Gamma_m(n_E/2) \Gamma_{n_H}\{(n_H + m + 1)/2\}\Gamma(1/2)}\mathrm{etr}\biggl(-\frac{1}{2}\Omega\biggl)\sum_{k=0}^{\infty}\frac{\{(n_H+n_E)/2\}_k}{(n_H/2)_k (m/2)_k}\frac{\Gamma(k+1/2)(\theta_1/2)^{k}}{k!}  \\
&\times \sum_{t=0}^{\infty} \sum_{\tau\in P^t_{n_H}}\sum_{\delta\in P^{k+t}_{n_H}}g_{(k),\tau}^{\delta}\frac{\{(m+1-n_E)/2\}_\tau(m/2)_{\delta}}{\{(n_H+m+1)/2\}_\delta t!}C_\delta (I_{n_H})\cdot x^{mn_H/2+k+t}
\end{align}
\begin{proof}
From $\mathrm{rank}(\Omega)=1$, we have $C_\kappa(\Omega/2)=C_{(k)}(\Omega/2)=(\theta_1/2)^k$.
Furthermore, using \eqref{zonal-identity}, we obtain the desired result by expressing the zonal polynomial $C_\kappa(I_m)$ for $\kappa=(k)$ as $C_\kappa(I_m)=2^kk!(m/2)_k/(2k)!=(m/2)_k\Gamma(1/2)/\Gamma(k+1/2)$.
\end{proof}
\end{corollary}
In the case of the linear alternative, the distribution \eqref{exact-rank1} is more useful than \eqref{prob-ell1} in terms of numerical computation, as the zonal polynomials $C_{(k)}(I_{n_H})$ and coefficients $g_{(k),\tau}^{\delta}$ in \eqref{exact-rank1} can be easily computed from \eqref{zonal-identity} and \eqref{g-coeffcient}, respectively.
Under the planar alternative, the zonal polynomials $C_\kappa(\Omega/2)$ in \eqref{prob-ell1} can be replaced by Legendre polynomials according to the bivariate property of zonal polynomials given in James~(1968).
If $\Omega=O$, we find the coefficients $g^\delta_{0, \tau}=1$ and $\delta=\tau$.
Then, the distribution \eqref{prob-ell1} under $H_0$ is expressed in the following form:
\begin{align}
\label{dist-H0}
\mathrm{Pr}(\ell_1 < x) &= \frac{\Gamma_{n_H}\{(n_H+1)/2\} \Gamma_m\{(n_H+n_E)/2\}}{\Gamma_m(n_E/2) \Gamma_{n_H}\{(n_H + m + 1)/2\}}
x^{mn_H/2}{_2F_1}\biggl(\frac{m+1-n_E}{2},\frac{m}{2};\frac{n_H+m+1}{2};xI_{n_H}\biggl),
\end{align}
where ${_2F_1}$ is the Gauss hypergeometric function of matrix arguments. 
Considering the transformation $q_1=\ell_1/(1-\ell_1)$, we obtain the exact distribution of the largest eigenvalue for the singular noncentral $F$ matrix defined as $F=B^{-1/2}AB^{-1/2}$ by 
\begin{align*}
\mathrm{Pr}(q_1 < x) =& \frac{\Gamma_{n_H}\{(n_H+1)/2\} \Gamma_m\{(n_H+n_E)/2\}}{\Gamma_m(n_E/2) \Gamma_{n_H}\{(n_H + m + 1)/2\}}\mathrm{etr}\biggl(-\frac{1}{2}\Omega\biggl)\sum_{k=0}^{\infty}\sum_{\kappa\in P^k_{n_H}}\frac{\{(n_H+n_E)/2\}_\kappa}{(n_H/2)_\kappa}\frac{C_\kappa(\Omega/2)}{C_\kappa(I_m)k!}  \\
&\times \sum_{t=0}^{\infty} \sum_{\tau\in P^t_{n_H}}\sum_{\delta\in P^{k+t}_{n_H}}g_{\kappa,\tau}^{\delta}\frac{\{(m+1-n_E)/2\}_\tau(m/2)_{\delta}}{\{(n_H+m+1)/2\}_\delta t!}C_\delta (I_{n_H})\cdot (x/(1+x))^{mn_H/2+k+t}
\end{align*}
The above distribution under $H_0$ was provided in Shimizu and Hashiguchi~(2022a).

Next, we consider the exact distribution of the Bartlett-Nanda-Pillai trace statistic $V^{(n_H)}$.
Khatri and Pillai~(1968) derived the density of $V^{(m)}$ in terms of the product of zonal polynomials. 
The density function of the Bartlett-Nanda-Pillai trace for $m>n_H$ can also be obtained from the density function of $V^{(m)}$ by replacing $m\to n_H$, $n_H\to m$ and $n_E\to n_H+n_E-m$; 
\begin{align*}
f(x)=& \frac{\Gamma_{n_H}\{(n_H+n_E)/2\}}{\Gamma_{n_H}\{(n_H+n_E-m)/2\}\Gamma(n_Hm/2)}~\mathrm{etr}\biggl(-\frac{1}{2}\Omega\biggl)\sum_{k=0}^{\infty}\sum_{\kappa \in P_{n_H}^{k}}
\frac{\{(n_H + n_E)/2\}_\kappa }{(m/2)_\kappa}\frac{C_\kappa(\Omega/2)}{C_\kappa(I_{n_H}) k!}\\
&\times \sum_{t=0}^{\infty}\sum_{\tau\in P_{n_H}^{t}}\sum_{\delta \in P_{n_H}^{k+t}}  g_{\kappa,\tau}^{\delta} \frac{\{(m+1-n_E)/2\}_\tau (m/2)_\delta C_\delta(I_{n_H})}{(mn_H/2)_{k+t} t!} x^{mn_H/2+k+t-1} 
\end{align*}
Integrating $f(x)$ with respect to $x$, we have the distribution function of $V^{(n_H)}$ as
\begin{align}
\label{dist-Pillai}
\nonumber
P(V^{(n_H)} < x) =& \frac{\Gamma_{n_H}\{(n_H+n_E)/2\}}{\Gamma_{n_H}\{(n_H+n_E-m)/2\}\Gamma(n_Hm/2)}~\mathrm{etr}\biggl(-\frac{1}{2}\Omega\biggl)\sum_{k=0}^{\infty}\sum_{\kappa \in P_{n_H}^{k}}
\frac{\{(n_H + n_E)/2\}_\kappa }{(m/2)_\kappa}\frac{C_\kappa(\Omega/2)}{C_\kappa(I_{n_H}) k!}\\
&\times \sum_{t=0}^{\infty}\sum_{\tau\in P_{n_H}^{t}}\sum_{\delta \in P_{n_H}^{k+t}}  g_{\kappa,\tau}^{\delta} \frac{\{(m+1-n_E)/2\}_\tau (m/2)_\delta C_\delta(I_{n_H})}{(n_H m/2+k+t)(mn_H/2)_{k+t} t!} x^{mn_H/2 +k+t}.
\end{align}
We note that \eqref{dist-Pillai} is convergent for $0<V^{(n_H)}<1$.
\begin{corollary}
Let $H\sim W_m(n_H,\Sigma,\Omega)$ and $E\sim W_m(n_E, \Sigma)$, where $m>n_H$, $H$, and $E$ are independent;
If $\mathrm{rank}(\Omega)=1$ and $\theta_1$ is the eigenvalue of $\Omega$, the exact distribution function of $V^{(n_H)}$ under the linear alternative is given as 
\begin{align}
\label{dist-V}
\nonumber
P(V^{(n_H)} < x) =& \frac{\Gamma_{n_H}\{(n_H+n_E)/2\}}{\Gamma_{n_H}\{(n_H+n_E-m)/2\}\Gamma(n_Hm/2)\Gamma(1/2)}~\mathrm{etr}\biggl(-\frac{1}{2}\Omega\biggl)\sum_{k=0}^{\infty}
\frac{\{(n_H + n_E)/2\}_k }{(m/2)_k (n_H/2)_k}\frac{\Gamma(k+1/2)(\theta_1/2)^k}{k!}\\
&\times \sum_{t=0}^{\infty}\sum_{\tau\in P_{n_H}^{t}}\sum_{\delta \in P_{n_H}^{k+t}}  g_{(k),\tau}^{\delta} \frac{\{(m+1-n_E)/2\}_\tau (m/2)_\delta C_\delta(I_{n_H})}{(n_H m/2+k+t)(mn_H/2)_{k+t} t!} x^{\frac{1}{2}n_H m+k+t}.
\end{align}
\begin{proof}
This is similar to the proof of Corollary~\ref{exact-rank1kei}.
\end{proof}
\end{corollary}

\section{Numerical computation}
\label{sec:04}
\subsection{Computation of a percentile point and time}
In this section, we calculate the derived distribution and obtain the exact powers of Roy's largest root and Bartlett-Nanda-Pillai trace statistics.
If $r=(m-{n_E}-1)/2$ is a positive integer, the distribution \eqref{prob-ell1} is a finite summation of $t$.
The truncated distribution of \eqref{prob-ell1} up to the $K$th degree with $r=(m-{n_E}-1)/2$ is represented by 
\begin{align}
\label{probfinite-ell1}
\nonumber
F_K(x) =& \frac{\Gamma_{n_H}\{(n_H+1)/2\} \Gamma_m\{(n_H+n_E)/2\}}{\Gamma_m(n_E/2) \Gamma_{n_H}\{(n_H + m + 1)/2\}}\mathrm{etr}\biggl(-\frac{1}{2}\Omega\biggl)\sum_{k=0}^{K}\sum_{\kappa\in P^k_{n_H}}\frac{\{(n_H+n_E)/2\}_\kappa}{(n_H/2)_\kappa}\frac{C_\kappa(\Omega/2)}{C_\kappa(I_m)k!}  \\
&\times \sum_{t=0}^{rn_H} \sum_{\tau\in P^t_{n_H}}\sum_{\delta\in P^{k+t}_{n_H}}g_{\kappa,\tau}^{\delta}\frac{\{(m+1-n_E)/2\}_\tau(m/2)_{\delta}}{\{(n_H+m+1)/2\}_\delta t!}C_\delta (I_{n_H})\cdot x^{mn_H/2+k+t}
\end{align}
The empirical distribution based on a $10^6$ Monte Carlo simulation is denoted by $F_{\mathrm{sim}}$.
Table~\ref{Percentilepoints}-(a) and (b) shows the percentile points between $F_K$ and $F_{\mathrm{sim}}$ in the linear and planer cases, respectively. 
We use Algorithm~3 in Shimizu and Hashiguchi~(2022b) to calculate $g^\delta_{\kappa, \tau}$ in \eqref{probfinite-ell1}.
Table~\ref{Percentilepoints} indicates that $F_{K}$ has three-decimal-place precision and achieves the desired accuracy. 
\begin{table}[H]
\caption{Percentile points of the truncated distribution of \eqref{prob-ell1}} 
\label{Percentilepoints}
\begin{center}
\begin{tabular}{c}
    \begin{minipage}[c]{0.4\hsize}
      \begin{center}
       \captionsetup{labelformat=empty,labelsep=none}
         \subcaption{$m=6, {n_H}=2, {n_E}=3, \theta_1=9$}
{\begin{tabular}{@{}cccc@{}} \toprule
$\alpha$&${{F^{-1}_\mathrm{sim}}}(\alpha)$ &$F^{-1}_{12}(\alpha)$  \\  \toprule
0.05	 &0.400& 0.400\\
0.10	 &0.454& 0.454\\
0.50	 &0.641& 0.641\\
0.90	 &0.799& 0.799\\ 
0.95	 &0.835& 0.835\\
\noalign{\smallskip}\hline
\end{tabular}}
 \end{center}
  \end{minipage}
  \begin{minipage}[c]{0.4\hsize}
          \begin{center}
        \captionsetup{labelformat=empty,labelsep=none}
           \subcaption{$m=6, {n_H}=2, {n_E}=3, \theta_1=9, \theta_2=3$}
{\begin{tabular}{@{}cccc@{}} \toprule
$\alpha$&${F_\mathrm{sim}^{-1}}(\alpha)$ &$F_{16}^{-1}(\alpha)$  \\  \toprule
0.05	 &0.435& 0.435\\
0.10	 &0.487& 0.487\\
0.50	 &0.664&0.664\\
0.90	 &0.811&0.811\\ 
0.95	 &0.844&0.844\\ 
\noalign{\smallskip}\hline
\end{tabular}}
        \end{center}
    \end{minipage}
  \end {tabular}
  \end{center}
\end{table}
If $r=(m-{n_E}-1)/2$ is a positive integer, the distribution \eqref{dist-H0} for $m=6, {n_H}=2, {n_E}=3$, and $\Omega=O$ can be represented by a finite series as follows:
\begin{align}
\label{finite}
f(x)=& 715 \left(\frac{x^{14}}{13}-\frac{112 x^{13}}{143}+\frac{1526 x^{12}}{429}-\frac{368
   x^{11}}{39}+\frac{34378 x^{10}}{2145}-\frac{7664 x^9}{429}+\frac{421
   x^8}{33}-\frac{16 x^7}{3}+x^6\right).
   \end{align}
 
   The $95$th percentile point of \eqref{finite} is 0.737.
  Thus, we obtain the exact powers as 0.231 and 0.276 for cases (a) and (b) in Table~\ref{Percentilepoints}, respectively.

Next, we consider the computation of \eqref{exact-rank1}. 
It is difficult to compute \eqref{prob-ell1} for large $p$ due to the complexity of computing zonal polynomials.
Using \eqref{zonal-identity} and \eqref{g-coeffcient}, the computation of \eqref{exact-rank1} can be performed more easily than \eqref{prob-ell1}. 
 We evaluate the computation time of \eqref{exact-rank1} by considering the truncated distribution as well as \eqref{probfinite-ell1}.
Table~\ref{comp-time} shows the computation time for the truncated distribution of \eqref{exact-rank1}, fixed $\theta_1=10$, and $K=20$.
The symbol ``-"  indicates a calculation time of more than one hour.
We find that the computation is easy to perform when the numbers of terms $rn_H$ and groups $p$ are small.
It takes less than two seconds to compute for \eqref{exact-rank1} but more than 100 seconds for \eqref{probfinite-ell1} with parameters $m=3$, $n_i=9$, and $p=3$.
\begin{table}[H]
  \begin{center}
\caption{Computation time for \eqref{prob-ell1} and \eqref{exact-rank1} for $\theta_1=10$ and $K=20$} 
\label{comp-time}
{\begin{tabular*}{13.6cm}{@{}cccccccc@{}} \toprule
$\begin{tabular}{c} Dim. \\ $m$ \end{tabular} $&
$\begin{tabular}{c} Samples\\ $n_i$ \end{tabular} $&
$\begin{tabular}{c} Groups\\ $p$ \end{tabular} $&
$\begin{tabular}{c} Terms \\$rn_H$ \end{tabular} $&
$\begin{tabular}{c} Comp. time of \eqref{exact-rank1}\\ (in seconds)\end{tabular} $&
$\begin{tabular}{c} Comp. time of \eqref{prob-ell1} \\ (in seconds)\end{tabular} $&\\ \toprule
3&5&3 &8& 0.191&11.63\\
3&9&3 &20&1.575&108.6 \\
5&4&4&6& 0.431&-	\\ 
5&6&4 &21& 6.004&-\\
7&3&5&4&0.088 	&-\\ 
7&5&5 &28& 20.07&-\\
9&3&6&5& 	0.140&- \\
9&4&6 &20& 12.76&-	\\
\noalign{\smallskip}\hline
\end{tabular*}}
 \end{center}
\end{table}
\subsection{Power caluculation}
\begin{figure}[H]
\begin{center}
\includegraphics[width=8cm]{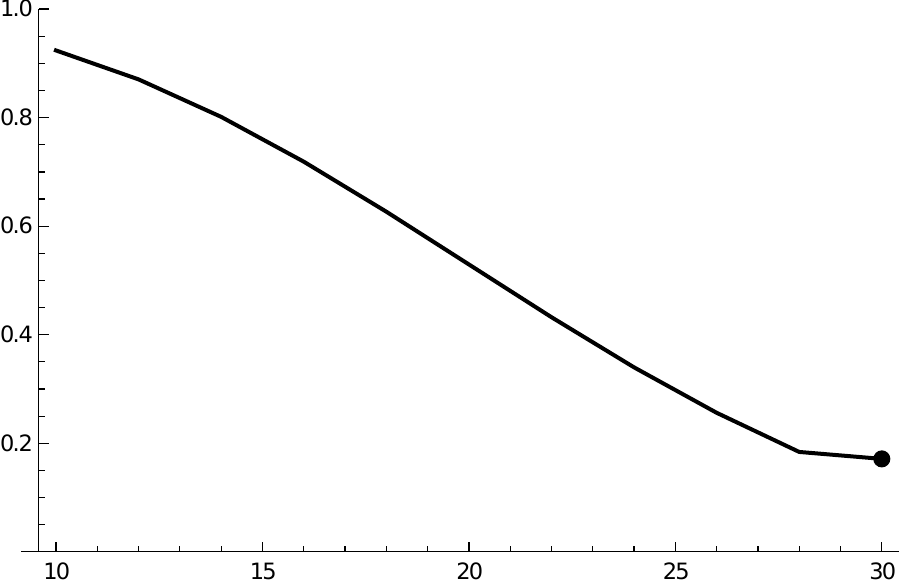}
\rlap{\raisebox{34ex}{\kern-24.em{\large $\mathrm{Power}$}}}%
\rlap{\raisebox{.15cm}{\kern0cm{\large $m$}}}
\caption{Power of Roy's test in the linear case with $\theta_1=40, n_i=12, p=3$}
 \label {fig1}
\end{center}
\end{figure}

We compare the exact powers of \eqref{exact-rank1} to the approximate powers in the linear case. 
The approximate powers can be obtained from Johnstone and Nadler's~(2017) results.
Table~3 presents a comparison between exact and approximate powers for the non-central parameter $\theta=10, 20, 40$. 
We confirm that Johnstone and Nadler's approach is not accurate when $\theta_1$ is small and $m$ is large. 
Therefore, our result is useful for obtaining exact powers for any dimension and noncentral parameter. 
\begin{table}[H]
  \begin{center}
\caption{Power comparison in the linear case} 
{\begin{tabular}{@{}cccccc@{}} \toprule
$\begin{tabular}{c} Dim. \\ $m$ \end{tabular} $&
$\begin{tabular}{c} Samples \\ $n_i$ \end{tabular} $&
$\begin{tabular}{c} Groups\\ $p$ \end{tabular} $&
$\begin{tabular}{c} Noncentrality \\ $\theta_1$ \end{tabular} $&
Exact&
\begin{tabular}{c} Johnstone and \\ Nadler's value \end{tabular}   \\  \toprule
3&3&5	 &10&0.229&0.182\\
3&3&5	 &20&0.465&0.437\\
3&3&5	 &40&0.811&0.812\\
7&6&4&10& 	0.210&0.139\\ 
7&6&4 &20& 0.453&0.407\\
7&6&4 &40&0.830&0.837\\
14&12&3&10&0.192&0.123\\ 
14&12&3&20& 0.415&0.374\\
14&12&3&40&0.801&0.823\\
\noalign{\smallskip}\hline
\end{tabular}}
 \end{center}
\end{table}
Next, we compare the exact powers in the linear case for the statistics between Roy's largest root and the Bartlett-Nanda-Pillai trace.
Table~\ref{powersR-V} shows the exact powers of Roys's test and Bartlett-Nanda-Pillai's test by \eqref{exact-rank1} and \eqref{dist-V} for a fixed $n_i=10~(i=1, 2, 3)$ and $p=3$. 
We find that the powers of Roy's test do not exceed those of the Bartlett-Nanda-Pillai trace test when the parameter $\theta_1$ is small.
\begin{table}[H]
  \begin{center}
\caption{Power comparison of two test statistics in the linear case with $n_i=10, p=3$}
\label{powersR-V}
{\begin{tabular*}{6.6cm}{@{}cccc@{}} \toprule
$\begin{tabular}{c} Dim. \\ $m$ \end{tabular} $&
$\begin{tabular}{c} Noncentrality \\ $\theta_1$ \end{tabular} $&
Roy&
Pillai ~\\  \toprule
4&3&0.147&0.148~~\\
6&3&0.118&0.120~~\\
8&3&0.101&0.105~~\\
4&10&0.475&0.442~~\\
6&10&0.396&0.363~~\\
8&10&0.292&0.267~~\\
\noalign{\smallskip}\hline
\end{tabular*}}
 \end{center}
\end{table}
Finally, we examine the exact powers under the linear alternative with an increase in dimensionality.
Fig~\ref{fig1} illustrates the exact powers of Roy's test for a fixed $\theta_1=40$, $n_i=12~(i=1, 2, 3)$, and $p=3$. 
Referring to Fig~\ref{fig1}, the exact powers decrease as $m$ increases; however, Roy's test is more powerful than the Bartlett-Nanda-Pillai trace test, especially when non-centrality is large in the linear alternative.
A power comparison for large dimensions with well-known tests under other alternatives should be conducted in future studies.

\section*{Acknowlegments}
The first author has been partially supported by Grant-in-Aid for JSPS Fellows (No.22J11185).




\end{document}